\documentclass{article}  

\title{An Algorithm for Constructing Cohomological Series Solutions
of Holonomic Systems}
\author{Nobuki Takayama}
\date{September 23, 2003, Revised January 15, 2004}

\def\vvv#1{  }

\def\pd#1{ \partial_{#1} }

\def\DD{ {\cal D} }

\def\OOf{ {\hat {\cal O}} }
\def\qed{ Q.E.D.  \par \bigbreak }
\def\xx{ {\bf x} }
\def\yy{ {\bf y} }
\def\Ext{ {\cal E}xt}

\newtheorem{theorem}{Theorem}[section]
\newtheorem{algorithm}[theorem]{Algorithm}
\newtheorem{example}[theorem]{Example}

\begin{document}
\maketitle

\section{Introduction}

Let $\DD$ be the sheaf of analytic differential operators of $n$ variables
$x_1, \ldots, x_n$.
We consider a left ideal $I$ of $\DD$ generated by 
$\{ \ell_1, \ldots, \ell_m \}$ which are in the Weyl algebra 
$D={\bf C}\langle x_1, \ldots, x_n, \pd{1}, \ldots, \pd{n} \rangle$.
If no confusion arises, we also denote by $I$ the left ideal 
$D \cdot \{ \ell_1, \ldots, \ell_n \}$ in $D$.
Assume that $\DD/I$ is holonomic.
It was proved by M.Kashiwara that the germs of the $k$-th extension group
$\Ext_\DD^k(\DD/I, \OOf)$
form a finite dimensional vector space over the field of complex numbers ${\bf C}$
\cite{kashiwara}.
We note that the vector space is called a $k$-th order 
(cohomological) solution space.
In \cite{oaku-takayama}, an algorithm by which to determine the dimension 
of the vector space was given.
In the present paper, 
we will present an algorithm by which to construct a basis of this vector space
in a free module over the formal power series.
In \cite{oaku-takayama}, we studied
the adapted free resolution of $D/I$
and an algorithm of computing restrictions of $D$-modules.
The algorithm of evaluating the dimension of the germ of the $k$-th
extension group was an immediate application of 
Cauchy-Kowalevski-Kashiwara's theorem
on the restriction of the $\DD$-module $\DD/I$ to the origin
and the $k$-th extension group.
In the present paper, we will explicitly construct matrix representations
of boundary operators of complexes appearing in a proof of the CKK Theorem
to construct series solutions.

Let
$$ \cdots \stackrel{\psi_{i+1}}{\rightarrow} D^{b_i}
         \stackrel{\psi_{i}}{\rightarrow}   D^{b_{i-1}} \rightarrow
  \cdots \rightarrow D \rightarrow  M \rightarrow 0 $$
be a free resolution of $M=D/I$.
Then, for a left $D$-module $N$, the vector space
$$ \frac{ {\rm Ker}\, \left( {\rm Hom}_D(D^{b_i},N) \rightarrow 
                             {\rm Hom}_D(D^{b_{i+1}},N) \right) }
        { {\rm Im}\,  \left( {\rm Hom}_D(D^{b_{i-1}},N) \rightarrow 
                             {\rm Hom}_D(D^{b_{i}},N)\right)  }
$$
is denoted by ${\rm Ext}_D^i(M,N)$ and
is called the $k$-th extension group.
When $M$ and $N$ are holonomic $D$-modules,
H.Tsai and U.Walther \cite{tsai-walther} presented an algorithm by which to
determine a basis of the extension group.
Let us consider
the sheaf of $k$-th extension group $\Ext_\DD^k(\DD/I,\OOf)$,
which can be defined in a similar way.
In this case, $\OOf$ is not holonomic over $D$,
and hence we cannot apply their algorithm; we need a different approach.

When $n=1$, it is relatively easy to determine bases of
$\Ext_\DD^k(\DD/I, \OOf)_0$, $k=0, 1$
where $I = \DD \cdot \ell$.
Consider the free resolution
$$ 0 \longrightarrow \DD \stackrel{\cdot \ell }{\longrightarrow}
                     \DD \stackrel{{\rm id}}{\longrightarrow} \DD/I
     \longrightarrow 0.
$$ 
By applying ${\rm {\cal H}om}_\DD(\cdot, \OOf)$-functor,
we have the complex 
$$ 0 \longleftarrow \OOf \stackrel{\ell \cdot}{\longleftarrow}
                    \OOf \longleftarrow 0.
$$
Hence, we have
$$ \Ext_\DD^0(\DD/I,\OOf) = \{ f \in \OOf \,|\, \ell \bullet f = 0 \}, \quad
   \Ext_\DD^1(\DD/I,\OOf) = \OOf/\ell \bullet \OOf
$$
following the definition of $\Ext$.
Algorithmic methods by which to determine bases of the vector spaces above
are well-known. 
Among these, we would like to examine in greater detail the method
explained in the introductory text book by T.Oaku
\cite{oaku-asakura}.
The key in this method is to regard $\OOf$ as an infinite
dimensional vector space over ${\bf C}$ and to
regard the operator $\ell$ as a linear map.
Oaku uses a $b$-function to reduce a problem of infinite dimensional vector spaces
into a problem of finite dimensional vector spaces.
Our method is a natural generalization of this method.

The motivation for the present study is the problem of constructing series solutions
in Gevrey classes of an ${\cal A}$-hypergeometric system
\cite{castro-takayama}. 
See also \cite{iwasaki} for the same problem for Lauricella hypergeometric
functions.
We hope that our algorithm can be applied to this problem.
In April 2003, the author discussed this problem with F.J.Castro-Jim\'enez, 
who presented an explicit computation of $k$-th order 
solutions of the hypergeometric system for $A=(1,2)$.
This example was exciting and has been the motivation for the present study.

\section{Orders in $\OOf$}

Let $K$ be a field and we consider the ring of the formal power series
$\OOf_0 = K[[x_1, \ldots, x_n]]$ 
in $n$-variables.
We will omit the subscript $0$ of $\OOf$ in the sequel.
When we apply the results of this section for $k$-th order solutions, 
$K$ is assumed to be ${\bf C}$.
We regard $\OOf$ as an infinite dimensional vector space
over $K$ by sorting monomials of $x$ by an order.
In other words, we identify $\OOf$ with a vector space
by sorting the coefficients of power series by the order.

Let us firstly consider when $n=1$.
If we sort the monomials as $1, x, x^2, x^3, \ldots, $,
then we have the canonical morphism $\sigma$ from $K^\infty$
to $\OOf$ as 
$$
\begin{array}{rcl}
   K^\infty := \{ c=(c_0,c_1, c_2, \ldots)\,|\, c_i \in K \} 
 & \stackrel{\sigma}{\longrightarrow} & K[[x]] \cr
  c
 & \stackrel{\sigma}{\longmapsto} &
 f = \sum \frac{c_k}{k!} x^k \cr
\end{array}
$$
In this example, the coefficients $c_k$ is sorted as
$c_0, c_1, c_2, \ldots, $.

When $n > 1$, we have several natural choices to sort monomials in $x$
depending on a weight vector and a term order;
let $w \in {\bf R}^n$ be a weight vector satisfying $w_i > 0$
and 
$\prec$ a term order in ${\bf Z}_{\geq 0}^n$.
We sort monomials in $x$ and consequently the coefficients $c_k$
of power series by the order $\prec_w$.
For example, when $w=(1, \ldots, 1)$  and $\prec$ is lexicographic order,
we sort the coefficients $c_k$ as 
$c=(c_{0 \cdots 0},c_{1 0 \cdots 0}, \ldots, c_{0 \cdots 0 1}, c_{2 0 \cdots 0}, \ldots )$.
In this case, we define the canonical morphism $\sigma$ by
$$
\begin{array}{rcl}
   K^{n\infty} := \{ c=(c_{0 \cdots 0}, c_{10 \cdots 0}, \ldots )  \,|\, 
       c_\alpha \in K \} 
 & \stackrel{\sigma}{\longrightarrow} & K[[x_1, \ldots, x_n]] \cr
  c
 & \stackrel{\sigma}{\longmapsto} &
 f = \sum \frac{c_\alpha}{\alpha!} x^\alpha \cr
\end{array}
$$
Here, $\alpha! = \alpha_1 ! \cdots \alpha_n!$,
$x^\alpha = x_1^{\alpha_1} \cdots x_n^{\alpha_n}$.

Finally, let us discuss how to encode 
$\OOf^r$ as an infinite dimensional vector space over $K$.
In this case, we use the degree shift vector ${\bf s} \in {\bf Z}^r$
in addition to a weight vector $w$ and a term order $\prec$.
Let ${\bf e}_i x^\alpha$ be the element of $\OOf^r$
where ${\bf e}_i$ is the $i$-th standard vector.
As in the theory of Gr\"obner basis for $D$ \cite{oaku-takayama},
we define
${\rm ord}_w[{\bf s}](e_i x^\alpha) = w \cdot \alpha + s_i$.
To define the canonical morphism $\sigma$,
we sort ${\bf e}_i x^\alpha$ lexicographically by 
${\rm ord}_w[{\bf s}]$, $-i$ and $\prec$ for $\alpha$.
The canonical morphism is defined as follows
\begin{equation} \label{eq:encode}
\begin{array}{rcl}
   \left(K^{n\infty}\right)^r = \{ c=(c_{i\alpha})  \,|\, 
    c_{i\alpha} \in K, \} 
 & \stackrel{\sigma}{\longrightarrow} & K[[x_1, \ldots, x_n]]^r \cr
  c
 & \stackrel{\sigma}{\longmapsto} &
 f = \sum \frac{c_{i\alpha}}{\alpha!} x^\alpha {\bf e}_i \cr
\end{array}
\end{equation}
Here,  $i=0, \ldots, r-1, \alpha \in {\bf N}_0^n$ 
and
${\bf e}_0, \ldots, {\bf e}_{r-1}$ are standard unit vectors.

We denote by 
$\left( K^{n\infty} \right)^r [{\bf s}]_{\leq m} $
the image by $\sigma^{-1}$ of ${\bf e_i} x^\alpha$ of which 
${\rm ord}_w[{\bf s}]$-degree is less than or equal to $m$.
There exist the following natural projections
$$
\begin{array}{crcl}
 \tau_m : & \left( K^{n\infty} \right)^r [{\bf s}]
          & \longrightarrow 
          & \left( K^{n\infty} \right)^r [{\bf s}]_{\leq m}  \cr
 \tau_{m'm} : & \left( K^{n\infty} \right)^r [{\bf s}]_{\leq m'}
          & \longrightarrow 
          & \left( K^{n\infty} \right)^r [{\bf s}]_{\leq m}  \cr
\end{array}
$$
Here, $m' \geq m$ and
$\left( K^{n\infty} \right)^r [{\bf s}]$ is the union of
$\left( K^{n\infty} \right)^r [{\bf s}]_{\leq m}$ for the  natural numbers $m$ 
in $\left( K^{n\infty} \right)^r$.

We call $\tau$ the truncation map.

\begin{example} \rm
Put $w=(1,1)$ and $s=(0,0,-1)$.
We denote  ${\bf e}_i x_1^j x_2^k$ by $[i,j,k]$. 
We consider $(K^{2\infty})^3[(0,0,-1)]$.
The elements of which degree is less than or equal to $1$ are
$$ [0,0,0],
          [1,0,0],[2,1,0],[2,0,1],  \mbox{(degree 0, 4 elements )}$$ 
$$ [0,1,0],[0,0,1],[1,1,0],[1,0,1],[2,2,0],[2,1,1],[2,0,2]  
   \mbox{(degree 1, 7 elements)}. $$
The projection $\tau_{1,0}$ is a map from 
$K^{11}$ to $K^4$.
\end{example}

\section{Adapted resolution and induced linear map} \label{sec:adapted}

We fix a weight vector $w$ in the sequel.
We consider a complex of left $D$-modules with degree shifts
${\bf m} \in {\bf Z}^p$,
${\bf m}' \in {\bf Z}^q$,
${\bf m}'' \in {\bf Z}^r$
\begin{equation}  \label{equation:sharp}
D^p[{\bf m}] \stackrel{A}{\rightarrow} D^q[{\bf m'}]
               \stackrel{B}{\rightarrow} D^r[{\bf m''}]
\end{equation}
(see \cite{oaku-takayama} on the degree shift).
We suppose that the following two conditions are satisfied. \\
(1) ${\rm ker}(B)$ is generated by rows $A_i$ of the matrix $A$.
(2) ${\rm in}_{(-w,w)}[{\bf m'}]({\rm ker}(B))$ is generated by
    ${\rm in}_{(-w,w)}[{\bf m'}](A_i)$.
When these two conditions are satisfied, the complex is called
{\it adapted} at the object $D^q[{\bf m}']$.
A free resolution is called adapted if it is adapted at every object
in the complex.
The notion is introduced in \cite{oaku-takayama} and
a free adapted resolution can be constructed by a Gr\"obner basis method
for a given left $D$ module $D/I$, a weight vector $w$, and a term order
$\prec$, which is used as a tie-breaker.
See \cite{oaku-takayama-minimal} for an efficient construction algorithm.

Let us consider the ring of formal power series
$\OOf = K[[x_1, \cdots, x_n]]$, 
and a complex
\begin{equation} \label{eq:sharp2}
  D^p[{\bf s}] \stackrel{\cdot A}{\longrightarrow}
  D^q[{\bf t}] \stackrel{\cdot B}{\longrightarrow}
  D^r[{\bf u}]
\end{equation}
which is adapted at the middle object $D^q[{\bf t}]$.
Note that we suppose that $D^p$, $D^q$, $D^r$ are sets of row vectors.
By applying 
${\rm Hom}_{\DD}( \cdot, \OOf)$
to the complex,
we have
\begin{equation} \label{eq:complexO}
  \OOf^p[{\bf s}] \stackrel{A \bullet}{\longleftarrow}
  \OOf^q[{\bf t}] \stackrel{B \bullet}{\longleftarrow}
  \OOf^r[{\bf u}]
\end{equation}
where $\OOf^i$, $(i=p, q, r)$ are regarded as sets of column vectors.
By the canonical isomorphism $\sigma$ (\ref{eq:encode}) induced 
by the degree shift, the weight vector and the tie-breaking 
term order $\prec$,
we have the following complex of linear vector spaces
\begin{equation} \label{eq:complexV}
\left(K^{n\infty}\right)^r[{\bf u}] 
\stackrel{\cdot {\bar B}^T}{\longrightarrow}
\left(K^{n\infty}\right)^q[{\bf t}] 
\stackrel{\cdot {\bar A}^T}{\longrightarrow}
\left(K^{n\infty}\right)^p[{\bf s}]. 
\end{equation} 
Here, 
${\bar B}^T$ and ${\bar A}^T$  are
block upper triangular matrices with the elements in $K$.
The blocks are partitioned by ${\rm ord}_w[\cdot]$-degree.
The two properties will be key ingredients of our algorithm.
Since it is block upper triangular,
we obtain the following complex of finite dimensional vector
space by truncating to the degree $m$
\begin{equation} \label{eq:complexVf}
\left(K^{n\infty}\right)^r[{\bf u}]_{\leq m} 
\stackrel{\cdot \tau_m({\bar B}^T)}{\longrightarrow}
\left(K^{n\infty}\right)^q[{\bf t}]_{\leq m}  
\stackrel{\cdot \tau_m({\bar A}^T)}{\longrightarrow}
\left(K^{n\infty}\right)^p[{\bf s}]_{\leq m}  
\end{equation}

Let us illustrate the matrix representations of the boundary operators
${\bar A}^T$ and ${\bar B}^T$ in terms of boundary operators
in (5.2) of \cite{oaku-takayama}, which are used to compute restrictions
of $D$-modules.
We consider the quotient ${\rm Ker}\, (A \bullet)/{\rm Im}\, (B \bullet)$
in (\ref{eq:complexO}).
First, we note that
$${\rm Ker}\, (A \bullet) = \left\{ {\bf f} \,\left|\,
    A \bullet {\bf f} = 0, \ 
    {\bf f} = \pmatrix{ f_1 \cr
                        \cdot \cr
                        \cdot \cr
                        \cdot \cr
                        f_q }, \   f_i \in \OOf \right. \right\}
$$
and
$$ {\rm Im}\, (B \bullet) = \left\{ B {\bf g} \, \left|\,
    {\bf g} = \pmatrix{ g_1 \cr
                        \cdot \cr
                        \cdot \cr
                        \cdot \cr
                        g_r }, \   g_i \in \OOf \right. \right\}.
$$
We denote by $A^j_\ell$ the $(j,\ell)$-th element of the matrix $A$.
Then the condition $A \bullet {\bf f} = 0$ is equivalent to
\begin{equation}  \label{eq:ker}
  \left( \pd{}^i \bullet ( \sum_\ell A^j_\ell \bullet f_\ell ) \right)
  (0) = 0, \quad
  \mbox{ for all } i \in {\bf N}_0^n \mbox{  and  } 
  j = 1, \ldots, p.
\end{equation}
Define $a^{ij}_{k\ell}$ by
$$ {\rm normallyOrdered}(\pd{}^i A^j_\ell)_{|_{x=0}} =
   \sum_{k} a^{ij}_{k\ell} \pd{}^k, \ 
  a^{ij}_{k\ell} \in K.
$$
Here, ${\rm normallyOrdered}(L)_{|_{x=0}}$ means that \\
(1) expand $L \in D$ into normally ordered expression as 
$\sum a_{\alpha k} x^\alpha \pd{}^k $ \\
and then \\
(2) replace all $x_i$ by $0$. \\
In terms of $a^{ij}_{k\ell}$, the condition (\ref{eq:ker}) is
equivalent to
\begin{equation}  \label{eq:ker2}
  \sum_{k,\ell} a^{ij}_{k\ell} f_\ell^{(k)}(0) = 0, \quad
\mbox{ for all } i \in {\bf N}_0^n, \ j=1, \ldots, p
\end{equation}
where $f_\ell^{(k)}$ is the $k=(k_1, \ldots, k_n)$-th derivative 
of $f_\ell \in \OOf$.
Note that $c_{i\alpha}$ in (\ref{eq:encode})
is equal to $f_i^{(\alpha)}(0)$.
Therefore, under the morphism $\sigma$,
${\rm Ker}\, (A \bullet)$ is nothing but the kernel of the
matrix defined by $(a^{ij}_{k\ell})$.

Next, we consider $B$.
Let $B^j_\ell$ be the $(j,\ell)$-th element of $B$.
Define $b^{ij}_{k\ell}$ by
${\rm normallyOrdered}(\pd{}^i B^j_\ell)_{|_{x=0}} = 
  \sum b^{ij}_{k\ell} \pd{}^k $.
Then the $i$-th coefficient of the series expansion of 
$\sum_\ell B^j_\ell \bullet g_\ell$ is expressed as
$$ \left( \frac{1}{i!} \pd{}^i \bullet
  \left(\sum_\ell B^i_\ell \bullet g_\ell\right) \right) (0)
 = \frac{1}{i!} \sum_{k, \ell} b^{ij}_{k\ell} g^{(k)}_\ell(0).
$$
Therefore, under the morphism $\sigma$,
${\rm Im}\, (B \bullet)$ is nothing but the image of the matrix
defined by $(b^{ij}_{k\ell})$.

The matrices $(a^{ij}_{k\ell})$ and $(b^{ij}_{k\ell})$ agree with
those appearing to compute restrictions of $D$-modules
\cite[Theorem 5.3]{oaku-takayama}.
Let us explain this fact.
The boundary operators $A$ and $B$ induce the following $K$-linear maps
$$ \Omega \otimes_D D^p[{\bf s}]
  \stackrel{\cdot {\bar A}}{\longrightarrow}
   \Omega \otimes_D D^q[{\bf t}]
  \stackrel{\cdot {\bar B}}{\longrightarrow}
   \Omega \otimes_D D^r[{\bf u}]
$$
where
$ \Omega = D/(x_1 D + \cdots + x_n D) $.
The cohomology groups of this complex are called {\it restrictions}.
See \cite{oaku-takayama} for details on computing restrictions.
Let us construct matrix representations of ${\bar A}$ and ${\bar B}$.
We denote by ${\bf e}^j$ the $j$-th standard vector.
The vector space $\Omega \otimes_D D^p[{\bf s}]$ is spanned by
 $\pd{}^i {\bf e}^j $, ($i \in {\bf N}_0^n$, $j =1, \ldots, p$),
which is sent by the linear map $\cdot {\bar A}$ to
${\rm normallyOrdered}(\pd{}^i \sum_\ell A^j_\ell {\bf e}^\ell)_{|_{x=0}}
 = \sum_{k,\ell} a^{ij}_{k\ell} \pd{}^k {\bf e}^\ell$.
By sorting $\pd{}^i {\bf e}^\ell$ by the same order to sort
the coefficients of power series
(we use the correspondence 
 $ \pd{}^i {\bf e}^\ell \leftrightarrow x^i {\bf e}_{\ell-1}$), 
we conclude that the matrix representation of ${\bar A}$ agrees with
the transpose of the ${\bar A}^T$ in (\ref{eq:complexV}).
The same assertion holds for ${\bar B}$.
This relation will be used in the proof of Theorem \ref{theorem:1}.

Let $M=D/I$ be a left holonomic $D$-module
and 
$$ \cdots \stackrel{\psi_{i+1}}{\rightarrow} D^{b_i}
         \stackrel{\psi_{i}}{\rightarrow}   D^{b_{i-1}}
         \stackrel{\psi_{i-1}}{\rightarrow} 
  \cdots
$$
an adapted free resolution associated to a weight vector $w$.
We assume that the complex (\ref{eq:sharp2}) is a part 
of this adapted resolution.
Let $k_1$ be the maximal integral root
of the $b$-function of $M$ associated to the weight $(-w,w)$
\cite{oaku-takayama}.

\begin{theorem}   \label{theorem:1}
The truncation map
$$ \tau_{m',m}\  :  \ 
  \frac{{\rm Ker}\, \tau_{m'}({\bar A}^T)}
       {{\rm Im}\, \tau_{m'}({\bar B}^T)}
  \longrightarrow
  \frac{{\rm Ker}\, \tau_{m}({\bar A}^T)}
       {{\rm Im}\, \tau_{m}({\bar B}^T)}
$$
is an isomorphism of vector spaces when
$m' > m \geq k_1$.
\end{theorem}

{\it Proof}.
The matrix $\tau_{m'}({\bar A}^T)$ is block triangular.
Define submatrices $A_{ij}^T$ by
$$\tau_{m'}({\bar A}^T) = \pmatrix{A_{11}^T & A_{21}^T \cr
                                  {\bf 0}  & A_{22}^T \cr}, \quad
  \tau_m({\bar A}^T) = A_{11}^T.
$$
Submatrices $B_{ij}^T$ are defined analogously.
Since (\ref{eq:complexVf}) is a complex for any $m$,
we have
$ B_{11}^T A_{11}^T = 0$,
$ B_{11}^T A_{21}^T + B_{21}^T A_{22}^T = 0$,
$ B_{22}^T A_{22}^T = 0$.
We note that the boundary operators $A_{22}^T$ and $B_{22}^T$
agree with those in the complex (5.3) of \cite{oaku-takayama}.
Since the complex (5.3) in \cite{oaku-takayama} is exact when
$m' > m \geq k_1$, we have
${\rm Im}\, B_{22}^T = {\rm Ker}\, A_{22}^T$.
Let $(p,q)$ be in ${\rm Ker}\, \tau_{m'}({\bar A}^T)$;
we assume that
$ p A_{11}^T = 0$ and
$ p A_{21}^T + q A_{22}^T = 0$.
Then, we can define a natural projection
$$ \psi \ :\ {\rm Ker}\, \tau_{m'}({\bar A}^T) \ni (p,q)
  \mapsto p \in {\rm Ker}\, A_{11}^T.
$$
By utilizing the properties of the matrices 
$A_{ij}^T$, $B_{ij}^T$ stated above,
it is easy to check that $\psi$ induces a well-defined map ${\bar \psi}$
from  
${\rm Ker}\, \tau_{m'}({\bar A}^T)/{\rm Im}\, \tau_{m'}({\bar B}^T)$
to ${\rm Ker}\, \tau_m(A_{11}^T) / {\rm Im}\, \tau_m(B_{11}^T)$
and that ${\bar \psi}$ is injective.
It follows from Theorem 5.3 of \cite{oaku-takayama} that
the dimensions of
$  \frac{{\rm Ker}\, \tau_{m'}({\bar A}^T)}
        {{\rm Im}\, \tau_{m'}({\bar B}^T)}$
and
$  \frac{{\rm Ker}\, \tau_{m}({\bar A}^T)}
        {{\rm Im}\, \tau_{m}({\bar B}^T)}$
agree.
Therefore, we conclude that ${\bar \psi} = \tau_{m',m}$  
is an isomorphism.
\qed

The next theorem immediately follows from Theorem \ref{theorem:1}.
\begin{theorem} \label{theorem:2}
Retain the notation of the proof of Theorem \ref{theorem:1}.
When $m' > m \geq k_1$,
for 
 $c \in {\rm Ker}\, \tau_{m}({\bar A}^T)$,
there exists a vector $c'$ such that
 $(c,c') \in {\rm Ker}\, \tau_{m'}({\bar A}^T)$.
In other words,
the linear inhomogeneous equation with unknown vector $c'$
$$  c \cdot A_{21}^T + c' \cdot A_{22}^T = {\bf 0} $$
is always solvable for  $c$ satisfying
$c \cdot A_{11}^T = {\bf 0}$.
\end{theorem}

\vvv{
{\it Proof}.
It follows from Theorem \ref{theorem:1} that
$$(c,q) \mapsto c 
   \in {\rm Ker}\, \tau_{m}({\bar A}^T)/{\rm Im}\, \tau_{m}({\bar B}^T)
$$
is surjective.
Hence, there exists $d$ 
such that
$(c+\tau_{m}({\bar B}^T)d, q') \in {\rm Ker}\, \tau_{m'}({\bar A}^T)$.
Subtracting
$\tau_{m'}({\bar B}^T)((d,0))$ from 
$(c+\tau_{m}({\bar B}^T)d, q')$,
we have
$(c, -d B_{12}^T +q')  \in {\rm Ker}\, \tau_{m'}({\bar A}^T)$.
}

\section{Algorithm}

Put ${\cal M} = \DD/I$.
We are ready to state our algorithm to construct a basis 
of the $d$-th order solutions.
It follows from Theorems \ref{theorem:1} and \ref{theorem:2}.
\begin{algorithm} \rm   \label{algorithm:1}
Construction of $d$-th cohomological solution 
${\rm {\cal E}xt}_\DD^d({\cal M}, \OOf)$\\
{\it Step 1.}  Construct an adapted resolution $(\psi_i)$ of $D/I$.
We assume that the resolution is written as (\ref{eq:sharp2})
at the degree $d$.  Note that we have assumed $A=\psi_{d+1}$, $B=\psi_d$.  \\
{\it Step 2.} Let $k_1$ be the maximal integral root of the $b$ function
of $M$ with respect to the weight vector $(-w,w)$. \\
{\it Step 3.}
Obtain a basis of 
$$ 
  \frac{{\rm Ker}\, \tau_{k_1}({\bar A}^T)}
       {{\rm Im}\, \tau_{k_1}({\bar B}^T)}
$$
as a $K$-vector space.
We denote the basis by ${\bf c}^{(1)}, \ldots, {\bf c}^{(e)}$. \\
{\it Step 4.}
By repeating to solve the linear equation in Theorem \ref{theorem:2},
extend the vector
${\bf c}^{(i)}$ to an infinite dimensional vector 
${\bf c}^{(i)}_\infty $ in
 $\left( K^{n\infty} \right)^q[{\bf t}]$ \\
{\it Step 5.}
Output 
$ \sigma({\bf c}^{(1)}_\infty), \ldots, \sigma({\bf c}^{(e)}_\infty)
\in K[[x_1, \ldots, x_n]]$ as a basis of the solutions.
\end{algorithm}

\noindent
{\it Remark}. \\
(1) In Step 3, our implementation chooses a basis
in 
$  {\rm Ker}\, \tau_{k_1}({\bar A}^T) \cap
   \left({\rm Im}\, \tau_{k_1}({\bar B}^T) \right)'
$
where  $V'$ denotes the orthogonal complement  of the vector space $V$ 
by the standard innerproduct. \\
(2) Note that we have to truncate the iteration of Step 4
when we execute this algorithm on a computer. \\
(3) Let $J$ be a submodule of $D^\ell$.
A basis of $\Ext_\DD^d(\DD^\ell/J, \OOf)$ can be constructed 
in an analogous way.
\bigbreak

The proof of the following theorem follows from 
discussions in Section \ref{sec:adapted}.

\begin{theorem}  \label{theorem:3}
The set of power series
$ \sigma({\bf c}^{(1)}_\infty), \ldots, \sigma({\bf c}^{(e)}_\infty) $
is a basis of 
$$ {\rm {\cal E}xt}_\DD^d({\cal M}, \OOf) =
  \frac{ \left\{ {\bf f} \in \OOf^q \,\left|\, A \pmatrix{f_1 \cr \cdot \cr 
                                                     \cdot \cr f_q } = 0 \right. \right\}}
  { \left\{ \left. B \pmatrix{g_1 \cr \cdot \cr 
                      \cdot \cr g_r } \,\right|\, g_i \in \OOf \right\}}
$$
as a $K$-vector space.
In particular, when $d=0$, 
${\rm {\cal E}xt}_\DD^d({\cal M}, \OOf)$ is a basis of the classical power
series solutions since the denominator is $0$ in the expression above.
\end{theorem}

\begin{example} \rm
Put
$$I = D\cdot \{ x \pd{x}-x (x \pd{x}+y \pd{y}+2) (x \pd{x}+3), \  
                y \pd{y}-y (x \pd{x}+y \pd{y}+2) (y \pd{y}+5) \}$$
and consider $M=D/I$.
The ideal $I$ annihilates the Appell function 
$F_2(2,3,5,1,1,x,y)$.
An adapted resolution of $M$ with respect to the weight vector
$(-1,-1,1,1)$
is as follows.
$$
   0 \longrightarrow D^2[(-1,0)] \stackrel{A}{\longrightarrow}
                     D^3[(0,0,-1)] \stackrel{B}{\longrightarrow}
                     D[(0)] \longrightarrow M \longrightarrow 0
$$
Here,
$$ B =
\pmatrix{
 \underline{x \pd{x}}-x^3 \pd{x}^2-x^2 y \pd{x} \pd{y}-6 x^2 \pd{x}-3 x y \pd{y}-6 x  \cr 
 \underline{y \pd{y}}-x y^2 \pd{x} \pd{y}-y^3 \pd{y}^2-5 x y \pd{x}-8 y^2 \pd{y}-10 y  \cr 
 b_{31} \cr
  } 
$$
$$ A = 
 \pmatrix{
  a_{11} &  a_{12}   & 
  \underline{-1}+x y \pd{y}+6 x-18 x y-117 x^2 y-135 x y^2  \cr 
  \underline{-y \pd{y}}+15 y+45 x y & \underline{x \pd{x}}-9 x-27 x y 
  & \underline{1}   \cr 
}
$$
where \\
$ b_{31} =  \underline{-x^3 y \pd{x}^2 \pd{y}+x^2 y^2 \pd{x}^2 \pd{y}-x^2 y^2 \pd{x} \pd{y}^2+x y^3 \pd{x} \pd{y}^2+5 x^2 y \pd{x}^2-7 x^2 y \pd{x} \pd{y}+9 x y^2 \pd{x} \pd{y}-3 x y^2 \pd{y}^2}+15 x^3 y \pd{x}^2+6 x^2 y^2 \pd{x} \pd{y}-9 x y^3 \pd{y}^2+45 x^4 y \pd{x}^2+45 x^3 y^2 \pd{x} \pd{y}-27 x^2 y^3 \pd{x} \pd{y}-27 x y^4 \pd{y}^2+270 x^3 y \pd{x}-135 x^2 y^2 \pd{x}+135 x^2 y^2 \pd{y}-216 x y^3 \pd{y} +270 x^2 y-270 x y^2 $ \\
and \\
$a_{11} =  \underline{x y^2 \pd{x} \pd{y}-x y^2 \pd{y}^2+y^3 \pd{y}^2+5 x y \pd{x}-7 x y \pd{y}+8 y^2 \pd{y}-5 y}+33 x y^2 \pd{y}+60 x y+162 x^2 y^2 \pd{y}+135 x y^3 \pd{y}+315 x^2 y-270 x y^2-2565 x^2 y^2-2025 x y^3-5265 x^3 y^2-6075 x^2 y^3 $, \\
$a_{12} =  \underline{-x^3 \pd{x}^2-3 x y \pd{y}}+3 x-18 x^2 y \pd{x}-9 x^2 y \pd{y}-54 x^2+27 x y-117 x^3 y \pd{x}-135 x^2 y^2 \pd{x}-27 x^2 y^2 \pd{y}-27 x^2 y+1053 x^3 y+1701 x^2 y^2+3159 x^3 y^2+3645 x^2 y^3 $.

Since the $b$-function is $s$, the number $k_1$ is equal to $0$.
The dimension of
${\rm {\cal E}xt}_\DD^0({\cal M},\OOf)$ is $1$ and 
$$ \sum_{m,n} \frac{(2)_{m+n} (3)_m (5)_n}{m! n!} x^m y^m $$
is a basis.

The dimension of ${\rm {\cal E}xt}_\DD^1({\cal M},\OOf) $ is
equal to $2$ at the origin.
Let us construct a basis of this vector space.
The induced linear maps on the space of truncated power series
$\tau_1({\bar A}^T)$, 
$\tau_1({\bar B}^T)$
are
$$ \tau_1({\bar A}^T) = \pmatrix{
\underline{ 0} & \underline{ -5}& \underline{ 0} & 0& 60& 0& 0& 15 \cr
\underline{ 3} & \underline{ 0}& \underline{ 0} & -108& 27& 0& -9& 0 \cr
\underline{ -1}& \underline{ 0} & \underline{ 0}& 12& 0& 0& 1& 0 \cr
\underline{ 0}& \underline{ -1} &\underline{ 0}& 0& 7& 0& 0& 1 \cr
0& 0& 0& 0& 0& 0& 0& 0 \cr
0& 0& 0& 0& -7& 6& 0& -1 \cr
0& 0& 0& 6& 0& 0& 1& 0 \cr
0& 0& 0& 0& 0& 0& 0& 0 \cr
0& 0& 0& -1& 0& 0& 0& 0 \cr
0& 0& 0& 0& -1& 0& 0& 0 \cr
0& 0& 0& 0& 0& -1& 0& 0 \cr
}
 $$
$$ \tau_1({\bar B}^T) = \pmatrix{
\underline{ 0}& \underline{ 0} & \underline{ 0}& \underline{ 0} & -6& 0& 0& -10& 0& 0& 0 \cr
0& 0& 0& 0& 1& 0& 0& 0& 0& 0& 0 \cr
0& 0& 0& 0& 0& 0& 0& 1& 0& 0& 0 \cr
}
 $$
Submatrices standing for the underlined elements are equal to
$\tau_0({\bar A}^T)$, $\tau_0({\bar B}^T)$.
A basis of the vector space
$\frac{ {\rm Ker}\, \tau_0({\bar A}^T)}
      { {\rm Im}\, \tau_0({\bar B}^T)}
$
is
$ (1,0,0,-5), \ (0,1,3,0)$.
Here the monomials ${\bf e}_i x^j y^k  $, which is encoded as $[i,j,k]$,
are sorted as
$[[0,0,0],  [1,0,0],[2,1,0],[2,0,1]]$.
Then the $0$-th approximation of a basis of series solutions is
$$ (1,0,-5y),  (0,1,3x) $$
This solution can be extended to the $1$-th approximation
by solving the linear inhomogeneous equation
in Theorem \ref{theorem:2} for $m=0$ and $m'=1$
as
$$(\underline{ 1}+ 10  y  , \underline{ 0} , 
   \underline{ -5  y} -45  y  x +60  y^{2}/2!  )$$
$$(\underline{ 0} , \underline{ 1}+6  x  , 
   \underline{ 3  x} -36  x^{2}/2!  +27 y  x )$$
The set of indices is sorted as
$$[[0,0,0],[1,0,0],[2,1,0],[2,0,1],[0,1,0],[0,0,1],[1,1,0],[1,0,1],[2,2,0],[2,1,1],[2,0,2]]$$

Repeating the procedure for
$\tau_2({\bar A}^T)$, $\tau_2({\bar B}^T)$, 
we obtain the following basis of second approximate solutions
\begin{eqnarray*}
(&&1+10\yy +\frac{1601775000} {10611803}\xx^{2}  +\frac{574824600} {10611803}  \yy  \xx +180  \yy^{2}   ,\   \\
 && \frac{574824600} {10611803}   \yy  \xx +\frac{278581950} {10611803} \yy^{2} ,\   \\
 && -5  \yy -45 \yy  \xx +60 \yy^{2}  -\frac{38669400} {10611803}  \yy  \xx^{2}  -\frac{39525750} {10611803}  \yy^{2}   \xx +7020  \yy^{3}  )
\end{eqnarray*} 
\begin{eqnarray*}
(&&\frac{280878030} {10611803}  \xx^{2}  +\frac{344820564} {10611803}  \yy  \xx   ,\ 
 1+6\xx +72  \xx^{2}  +\frac{344820564} {10611803}\yy  \xx +\frac{2064283056} {10611803} \yy^{2}    ,\  \\
&& 3 \xx -36  \xx^{2}  +27 \yy  \xx -2376 \xx^{3}  +\frac{15126426} {10611803}  \yy  \xx^{2}  +\frac{13920984} {10611803} \yy^{2}   \xx )
\end{eqnarray*}
Here, we set $\xx^i \yy^j = x^i y^j/(i! j!)$.
\end{example}

\end{document}